\documentstyle[amscd]{amsart}

\numberwithin{equation}{section}

\newtheorem{theorem}{Theorem}[section]

\newtheorem{lemma}[theorem]{Lemma}
\newtheorem{corollary}[theorem]{Corollary}
\newtheorem{proposition}[theorem]{Proposition}

\theoremstyle{remark}
\newtheorem{remark}[theorem]{Remark}

\newtheorem{definition}[theorem]{Definition}
\newtheorem{example}[theorem]{Example}

\newcommand{\ad}{\operatorname{ad}}

\newcommand{\End}{\operatorname{End}}
\newcommand{\Ker}{\operatorname{Ker}}

\newcommand{\g}{{\frak g}}
\newcommand{\cx}{{\Bbb C}}

\begin{document}

\title{Manifolds with an $SU(2)$-action on the tangent bundle}
\author{Roger Bielawski}
\address{Department of Mathematics\\
University of Glasgow\\
Glasgow G12 8QW\\ Scotland.}
\email{R.Bielawski@@maths.gla.ac.uk}

\thanks{Work supported by an EPSRC Advanced Research Fellowship}


\begin{abstract} We study manifolds arising as spaces of sections of complex manifolds fibering over ${\Bbb C}P^1$ with normal bundle of each section isomorphic to ${\cal O}(k)\otimes {\Bbb C}^n$.
\end{abstract}

\maketitle
\thispagestyle{empty}

Any hypercomplex manifold can be constructed as a
space of sections of a complex manifold $Z$ fibering over ${\Bbb C}P^1$. The normal bundle  of each section must be the sum of ${\cal O}(1)$'s and this suggests that interesting geometric structures can be obtained  if we replace ${\cal O}(1)$ with other line bundles. Such structures have been introduced by many authors \cite{Aki, AG, Bai, DM, DM2, Gi, Ma} and are variously known as conic, Grassman, paraconformal or ${\cal P}$-structures.
   \par
Spaces of sections with normal bundle ${\cal O}(n)\oplus {\cal O}(n)$ have been studied recently,  in detail, by Maciej Dunajski and Lionel Mason \cite{DM, DM2}, and much of the present paper can be viewed as translation of their work from spinor language. 
\par
We adopt the point of view that spaces of sections of complex
manifolds fibering over ${\Bbb C}P^1$ are {\em manifolds with a fibrewise linear
action of $SU(2)$ on the tangent bundle}. The integrability condition says that the ideal in $(\Omega^\ast M)^\cx$ generated by the highest weight $1$-forms is closed, for any Borel subgroup. We call such manifolds {\em generalised hypercomplex manifolds} (or {\em $k$-hypercomplex manifolds} if the representation of $SU(2)$ splits into copies of the $k$-th symmetric power of the standard representation). These manifolds seem to be interesting from several points of view, quite apart from ``intrinsic worth". On the one hand, they provide a natural setting for certain integrable systems, ``monopoles", as we discuss in sections \ref{bundles} and \ref{monopoles}. These comprise the Bogomolny hierarchy of Mason and Sparling \cite{MS} and the self-dual hierarchy. On the other hand, even if one is interested primarily in hypercomplex or hyperk\"ahler manifolds, one may wish to consider these generalised hypercomplex manifolds, since (for odd $k$) they are foliated by hypercomplex submanifolds (see section \ref{blow}; this observation goes back to Gindikin \cite{BG}, see also \cite{DM,DM2}). 
\par
The paper is organised is as follows: in the next section we recall some basic facts about foliations and distributions. In section \ref{definition} we define the generalised hypercomplex manifolds and their twistor spaces. Section \ref{Quillen} is devoted to an interpretation of a construction of D. Quillen \cite{Q}, which we need later on. In sections \ref{bundles} and \ref{monopoles} we discuss the Ward correspondence and monopoles on $k$-hypercomplex manifolds. In the following section we show that  a $k$-hypercomplex manifold $M$ has, analogously to hypercomplex manifolds,  an  $S^2$-worth of certain integrable structures. For odd $k$, $M$ is always a complex manifold, but the $S^2$-worth of complex structures are combined quite differently from hypercomplex structures. In section \ref{blow} we show that $M$ is foliated by $(k-2i)$-hypercomplex submanifolds, for any $i<k/2$. In particular, as mentioned above, for odd $k$ we obtain a foliation by hypercomplex submanifolds. In section  \ref{ext} we construct a "hypercomplex extension" of a $k$-hypercomplex  manifold $M$, i.e. a hypercomplex manifold $\tilde{M}$ fibering over $M$ with the property that a monopole on $M$ is equivalent to a solution of self-duality equations on   $\tilde{M}$. In the next section we discuss maps between generalised hypercomplex manifolds, and in section \ref{form} symplectic structures and symplectic quotients. In the last section we give some examples -   the $k$-hypercomplex manifolds arising as moduli spaces of solutions to analogues of Nahm's equations seem to be particularly interesting.
\par
Finally, we should ask whether we could replace the action of $SU(2)$ on $TM$ with other compact semisimple Lie groups? Thus we ask for manifolds with a fibrewise $G$-action on $TM$ such that the ideal in $(\Omega^\ast M)^\cx$ generated by the highest weight $1$-forms is closed for any Borel subgroup of $G^\cx$. An example of such a manifold is the group $G$ itself. Much of the theory goes through for such ``{\em Borel-Weil manifolds}", but there is dearth of examples. In any case, twistor considerations show that for $G\neq SU(2)$, the canonical linear connection on $M$ (see Remark 5.7) is  {\em flat} (not necessarily torsion-free).

{\bf Remark on notation:} If X is a complex manifold, then $TX, T^\ast X,
\Omega^p(X)$ all denote holomorphic objects, i.e. the $(1,0)$-tangent and
-cotangent bundle and the sheaf of holomorphic $p$-forms, respectively. If $E$
is a holomorphic bundle on $X$, we denote by $E$ also the sheaf of its sections.

\section{Foliations and involutive structures\label{involutive}}

Here, we recall some basic facts about foliations and involutive structures (i.e. complex distributions). The basic references are \cite{Mol} and \cite{Tre}.
\begin{definition} A {\em foliated manifold} is a manifold $M$ of dimension $n$ modelled on the fibration ${\Bbb R}^q\times {\Bbb R}^{n-q}\rightarrow {\Bbb R}^q$. This means that there is a smooth atlas $\{U_i,\phi_i\}$, $\phi_i:U_i\rightarrow {\Bbb R}^q\times {\Bbb R}^{n-q}$, such that the transition functions $\phi_i\circ\phi_j^{-1}:U_i\cap U_j\rightarrow {\Bbb R}^q\times {\Bbb R}^{n-q}$ are of the form 
\begin{equation} {\Bbb R}^q\times {\Bbb R}^{n-q}\ni (x,y)\longmapsto (h_{ij}(x),p_{ij}(x,y))\in {\Bbb R}^q\times {\Bbb R}^{n-q}.\label{trans} \end{equation}
\label{foliation}
\end{definition}
Equivalently, a foliation is given by an integrable distribution ${\cal P}\subset TM$. The integral submanifolds (leaves) of  ${\cal P}$ determine a partition $F$ of $M$, and, following \cite{Mol}, we avoid the constant reference to the foliated atlas through the expedient of denoting a foliated manifold by $(M,F)$. We shall also refer to $F$ as a foliation (of codimension $q$). 
\par
A foliation is {\em simple} if it is defined by a submersion. In particular, the foliated atlas of a foliated manifold consists of simple open sets.

Let now $P$ be a property of manifolds which is determined by a reduction of the pseudogroup of diffeomorphisms of ${\Bbb R}^m$ to a subpseudogroup $\Pi=\Pi_m$ (e.g. complex, affine etc.). We say that a foliation $F$ is a {\em $P$-foliation} if every leaf has property $P$, i.e. if for every $x\in {\Bbb R}^q$, the functions $p_{ij}(x,\cdot)$ in \eqref{trans} belong to $\Pi_{n-q}$. We shall say that a foliation $F$ is a {\em transversely $P$-foliation} if the functions $h_{ij}$ in \eqref{trans} belong to $\Pi_q$.
\par
We also need the definition of bundles and forms in this setting:
\begin{definition} Let $(M,F)$ be a foliated manifold and let $E$ be a vector (or principal) bundle on $M$. $E$ is said to be a {\em foliated bundle} if the restriction of $E$ to any simple open subset $U$ of $M$ is a pullback of a bundle on the local quotient manifold $\bar{U}$.\label{f-bundle}\end{definition} 
\begin{definition} An $r$-form $\omega$ on a foliated manifold $(M,F)$ is called {\em basic}, if it is locally a pullback of a form on a local quotient manifold. In other words $\omega$ can be represented locally using only the ``transverse" coordinates $x$ in \eqref{trans}.\label{basic}\end{definition}

We now turn to complex distributions.

\begin{definition} An {\em involutive} structure on a smooth manifold $M$ is an
involutive vector subbundle ${\cal V}$ of $T^\cx M$, i.e. ${\cal V}$ satisfies
$[{\cal V}, {\cal V}]\subset{\cal V}$. \end{definition}

The dual description gives us a subbundle $T^\prime$ of complex $1$-form
which vanish on ${\cal V}$. If ${\cal V}$ is involutive, then $T^\prime$ is
{\em closed}, i.e. for any local smooth section $\phi$ of $T^\prime$ there
are sections $\psi_1,\dots,\psi_m$ of $T^\prime$ and smooth differential
$1$-forms $\omega_1,\dots,\omega_m$ such that $$
d\phi=\psi_1\wedge\omega_1+\cdots \psi_m\wedge\psi_m.$$
\begin{definition} An  involutive structure on a smooth manifold $M$
given by the vector bundle ${\cal V}$ is called {\em elliptic} if $T^\cx
M={\cal V} +\bar{\cal V}$ (i.e. $T^\prime\cap
\bar{T^\prime}=0$)\end{definition}
The importance of elliptic structures follows from the following
generalisation of the Newlander-Nirenberg Theorem (cf. \cite{Tre}, Theorem
VI.7.1):
\begin{theorem} Let ${\cal V}\subset  T^\cx M$ be an elliptic involutive
structure.
Then  ${\cal V}$ is integrable, i.e. every point $m\in M$ has a
neighbourhood with local coordinates $z_1,\dots,z_r,t_1,\dots,t_k$, where
$z_i\in{\Bbb C}$ such that $$T^\prime=\text{span}\{dz_1,\dots,dz_m\}.$$
Equivalently
$$ {\cal V}=\text{span}\left\{\partial/\partial{\bar z}_j,
\partial/\partial t_k \right\}.\label{NN}$$\end{theorem}

Let ${\cal F}={\cal V}\cap \bar{\cal V}\cap TM$, where $\bar{\cal V}$ is the
complex conjugate of ${\cal V}$. Then ${\cal F}$ is a (real) integrable
distribution on $M$ and the above theorem says the corresponding foliation is a transversely holomorphic foliation. In particular, 
if the space of leaves  $M/{\cal F}$ is Hausdorff, then
it is a complex manifold.

\begin{remark} If $M$ and the involutive structure ${\cal V}$ are real-analytic, we can extend ${\cal V}$ to an involutive subbundle of $T^{1,0}M^\cx$, where $M^\cx$ is a complexification of $M$. Thus an involutive structure on $M$ corresponds to a holomorphic foliation of $M^\cx$.\end{remark}

\section{Generalised hypercomplex manifolds \label{definition}}

Almost complex manifolds can be thought of as manifolds whose tangent bundle
admits a fibrewise action of $U(1)$ such that each tangent space is the
direct sum of the standard $2$-dimensional representations of $U(1)$.
Similarly, an almost hypercomplex manifold is a manifold whose tangent
bundle admits a fibrewise action of $SU(2)$ such that each tangent space is
the direct sum of the standard $4$-dimensional representations of $SU(2)$. A
less well-known example is that of an $f$-structure \cite{Y,Raw}, which
amounts to giving a fibrewise $S^1$-action on $TM$ such that each tangent
space decomposes into standard or trivial $S^1$-representations.

\subsection{Generalised almost hypercomplex structures}

Let $M$ be a smooth manifold. A {\em generalised almost hypercomplex structure}
on $M$ is a smooth fibrewise action of $SU(2)$ on $TM$ such that each
tangent space is isomorphic to  $V\otimes {\Bbb R}^n$, where
$V$ is a fixed non-trivial irreducible representation of $SU(2)$ (on a real
vector space). The complexified representation $V^{\Bbb C}$ is then one or two
copies of the $k$-th symmetric power of the standard $2$-dimensional unitary
representation of $SU(2)$, and we shall also call $M$ an {\em almost
$k$-hypercomplex manifold}.
\par
An almost $k$-hypercomplex manifold has dimension $m(k+1)$, where $m$ is even if $k$ is odd. The structure group of such a manifold reduces to the centraliser of $SU(2)$ in $GL\bigl(m(k+1),{\Bbb R})$, i.e. to $GL(m,{\Bbb R})$ if $k$ is even and to $GL(m/2,{\Bbb H})$ if $k$ is odd. Let $E=E_M$ be the complex vector bundle on $M$ associated to the the standard representation of $GL(n/2,{\Bbb H})$ or $GL(m,{\Bbb R})$ on ${\Bbb C}^m$ (in the second case, the representation
is the complexification of the standard real representation). Let $H$ be the trivial bundle with fibre $S^k\cx^2$. We have then a canonical isomorphism:
\begin{equation}  TM^\cx\simeq E_M\otimes H.\label{E}\end{equation}
For even $k$ there is a corresponding splitting of the real tangent bundle, but we prefer to treat both even and odd $k$ uniformly. 
\begin{remark} The splitting as in \eqref{E} with trivial bundle $H$ is called in literature a {\em right-flat almost Grassman structure} \cite{Ma, AG}. Obviously, if we have such a splitting and the bundles $E,H$ are equipped with either quaternionic or real structures (depending on the parity of $\dim M$), we can define an almost $k$-hypercomplex structure on $M$. Thus there is no difference between generalised almost hypercomplex structures and right-flat almost Grassman structures. The difference comes in when we consider integrability conditions.\end{remark} 
\par
We now turn to constructing  almost $k$-hypercomplex manifolds.
One way of realising the irreducible representations of $SL(2,\cx)$ is as
sections of line bundles over $\cx P^1$. Similarly an irreducible
representation of $SU(2)$ on a real vector space can be realised as the space of
{\em real} sections of an irreducible $\sigma$-bundle on  $\cx P^1$ \cite{Q}.
Here $\sigma$-bundle means a holomorphic bundle equipped with an
anti-holomorphic involution $\tau$ covering the antipodal map $\sigma$ on   $\cx
P^1$. An irreducible $\sigma$-bundle is isomorphic to either ${\cal O}(2k)$ or
${\cal O}(2k+1)\oplus {\cal O}(2k+1)$. Therefore  each tangent space of a
generalised almost hypercomplex manifold can be realised as the space
of real sections of a  $\sigma$-bundle $E$ on
${\Bbb C}P^1$ and one way of obtaining
manifolds with a generalised almost hypercomplex structure is as the space of sections of a complex
manifold fibering over ${\Bbb C}P^1$. Indeed, we have
\begin{proposition} Let $Z$ be complex manifold fibering over ${\Bbb C}P^1$
and equipped with an antiholomorphic involution $\tau$ covering the
antipodal map $\sigma$ on ${\Bbb C}P^1$. Suppose that there exists a holomorphic
and
$\tau$-invariant section of $Z\rightarrow {\Bbb C}P^1$ whose normal bundle is
isomorphic to ${\cal O}(k)\otimes {\Bbb C}^n$, $k>0$. Then the space of such
sections is a manifold of dimension $n(k+1)$,
equipped with a canonical almost $k$-hypercomplex
structure.\label{twistor}\end{proposition}
\begin{pf} The fact that the space $M$ of real sections with normal bundle
${\cal O}(k)
\otimes {\Bbb C}^n$ is a manifold of dimension $n(k+1)$ follows from the
Kodaira deformation theory, given that $h^1({\cal O}(k))=0$. To show that there
is a canonical almost
$k$-hypercomplex structure on $M$, consider the vertical bundle $N\subset TZ$,
i.e. the
kernel of $d\pi:TZ\rightarrow T{\Bbb C}P^1$. For any $m\in M$, i.e. a section
$m:{\Bbb C}P^1\rightarrow Z$, the normal bundle to $m$ is just the restriction
of $N$
to $m$. The tangent space $T^{\Bbb C}_mM$ is identified with
$H^0(m,N_{|_m})$. Let $L={\cal O}(k)$. We have a canonical identification of
bundles on $Z$: $$
N\simeq \left(N\otimes \pi^*(L^\ast)\right)\otimes \pi^*(L)$$ which, given
the fact that $N_{|_m}\simeq L\otimes {\Bbb C}^n$ leads to the decomposition
$$T_mM^{\Bbb C}\simeq H^0(m, N\otimes \pi^*(L^\ast))\otimes H^0(m,
\pi^*(L)).$$ 
In other words the tangent bundle $T^{\Bbb C}M$ decomposes as the
tensor product  of two bundles on $M$, $H^0(\cdot, N\otimes \pi^*(L^\ast)$
and $H^0(\cdot, \pi^*(L))$. This second bundle is trivial and so the action
of $SU(2)$ on $H^0({\Bbb C}P^1,L)$ induces a canonical fibrewise action of $SU(2)$
on $T^{\Bbb
C}M$. This action restricts to the real tangent bundle as the existence of
the real structure $\tau$ implies that the representation of $SU(2)$ on
$H^0({\Bbb C}P^1,
L\otimes {\Bbb C}^n)$ is real.\end{pf}
\begin{example} The $3$-sphere carries a canonical $2$-hypercomplex structure, defined by identifying $TS^3\simeq SU(2)\times {\frak su}(2)$ and considering the adjoint action of $SU(2)$ on the second factor. The integrability can be checked directly, or we can view $S^3$ as the space of real sections of ${\Bbb P}\bigl( {\cal O}(1)\oplus {\cal O}(1)\bigr)\simeq \cx P^1\times \cx P^1$. \end{example}

\subsection{Integrability of GHC-structures}

\begin{definition}
We shall say that a generalised almost hypercomplex structure is integrable if
it arises
locally as in the previous proposition, i.e. $M$ can be realised (at least
locally) as the space of real sections of a complex manifold $Z$ fibering
over ${\Bbb C}P^1$. We shall call an integrable generalised almost hypercomplex
structure  (resp. integrable almost $k$-hypercomplex structure) a generalised
hypercomplex structure (resp. a  $k$-hypercomplex structure). We shall often
abbreviate the words ``generalised
hypercomplex" to GHC.\label{GHC}\end{definition}
Before proceeding, it will be useful to describe the map $$\text{({\em an
irreducible representation of} $SL(2,\cx)$)}\mapsto \text{({\em a line bundle on
$\cx P^1$})}.$$
If $H$ is such a representation, then $SL(2,\cx)$ acts irreducibly on $U^\ast$.
Fix  an isomorphism $\cx P^1\simeq SL(2,\cx)/ B$,
where $B$ is a Borel subgroup. For every Borel subgroup $B_q$ corresponding to a
point $q\in \cx P^1$, let $l_q$ be the line of highest weight vectors for
$B_q$. This determines a holomorphic line bundle $\tilde{L}$ on ${\Bbb C}P^1$.
In
fact, we obtain an imbedding ${\Bbb C}P^1\rightarrow {\Bbb P}(U^\ast)$ and
$\tilde{L}$
is the pullback of the tautological bundle on ${\Bbb P}(U^\ast)$. It follows
that $L=\tilde{L}^\ast$ is positive and $H\simeq H^0({\Bbb C}P^1,L)$ as
representations of $SL(2,\cx)$.

We now describe the condition on integrability of generalised almost
hypercomplex structures. Consider the
action of $SL(2,\cx)$ on the complexified cotangent bundle $(T^\ast M)^\cx$. For
for any point $q\in {\Bbb C}P^1$ which corresponds to a  Borel subgroup $B$ of
$SL(2,\cx)$, we define the following subbundles of $T^{\Bbb C}M$:
\begin{itemize}
\item ${\cal U}_q$ -  the subbundle  of $(T^\ast M)^\cx$ corresponding to the
highest
weight;
\item ${\cal K}_q$ - the subbundle of $T^\cx M$ annihilated by  ${\cal U}_q$;
\item ${\cal F}_q={\cal K}_q\cap \overline{{\cal K}}_q\cap TM$ - a
distribution on $M$.\end{itemize}
Thus ${\cal K}_q$ consists of all but the lowest weight tangent vectors.
 We have:
\begin{theorem} An almost $k$-hypercomplex structure on a manifold $M$ is
integrable
if and only if for every $q\in \cx P^1$ the corresponding
subbundle ${\cal K}_q$ of $T^{\Bbb C}M$ is involutive, i.e. $[{\cal K}_q,
{\cal K}_q]\subset {\cal K}_q$.\label{integrable}\end{theorem}
\begin{pf}
Suppose that the almost generalised hypercomplex structure is integrable, i.e.
$M$ is
given (locally) as the parameter space of $\tau$-invariant sections of a
holomorphic manifold $Z$ fibering over ${\Bbb C}P^1$. Then $M$ has the natural
complexification $M^{\Bbb C}$ defined as the parameter space of all sections
of $Z$ with the normal bundle ${\cal O}(k)\otimes {\Bbb C}^n$. For $q\in {\Bbb
C}P^1$ consider
the natural holomorphic map $p_q$ from $M^{\Bbb C}$ to the fiber $Z_q$  of
$Z$ over $q$ given by intersecting a section with the fibre. The map $dp_q$ is
given by highest-weight $1$-forms (this follows from the proof of
Proposition \ref{twistor}). The kernel of
$dp_q$ is a subbundle ${\cal W}$ of $T^{1,0}M^{\Bbb C}$. At points of $M$ we
have a canonical identification $T^{\Bbb C}_m M\simeq T^{1,0}_m M^{\Bbb C}$
under which ${\cal W}\simeq {\cal K}_q$.
\par
Conversely, suppose that the subbundle ${\cal K}_q$ is involutive for any
point $q$ of $\cx P^1$. We define a subbundle ${\cal V}$ of the complexified
tangent bundle
to $M\times {\Bbb C}P^1$ as follows:
\begin{equation} {\cal V}_{(m,q)}=({\cal K}_q)_m\oplus T_q^{0,1}{\Bbb C}P^1,
\label{Inv}
\end{equation}
where ${\Bbb C}P^1$ is equipped with the canonical complex structure. We
claim that this subbundle is involutive. Let us choose a point $q$ in ${\Bbb
C}P^1$
and let $v$ be a local section of ${\cal K}_q$. Let $u(\zeta)$ be a local
holomorphic section of $SL(2,\cx)\rightarrow \cx P^1$ (in a neighbourhood of
$q$). Then, because ${\cal K}$ is a homogeneous bundle, $u(\zeta)v$ is a
section of ${\cal K}_{u(\zeta)q}$ for each $\zeta$. The vector fields
$u(\zeta)v$ and $\partial/\partial \bar{\zeta_i}$, $\zeta_i$ being local
holomorphic coordinates on ${\Bbb C}P^1$, generate the bundle ${\cal V}$.
Therefore
${\cal V}$ is involutive.
 This involutive structure is
elliptic and hence  (replacing $M$ by an open subset, if Hausdorffness is a
problem) $(M\times {\Bbb C}P^1)/{\cal F}$, ${\cal F}={\cal V}\cap \bar{\cal
V}\cap TM$, is a complex manifold. This is our twistor space $Z$. It is
clear that $Z$ is a complex fibration over ${\Bbb C}P^1$ with the fiber at $q$
equal
to $M/{\cal F}_q$. Points of $m$ give rise to sections of $M\times
{\Bbb C}P^1\rightarrow {\Bbb C}P^1$, which then descend to sections of $Z$. That
these sections
have correct normal bundle follows from the remarks following Definition
\ref{GHC}.\end{pf}

\begin{remark} One can now generalise the notion of GHC-structures by allowing
group actions on $TM$ such that each tangent space splits into {\em
nonequivalent} irreducible representations and using
 Theorem \ref{integrable} as the definition. These are the ${\cal P}$-structures of Simon Gindikin \cite{Gi}.  \end{remark}

\subsection{The twistor foliation of a GHC-manifold}

The proof of the above theorem shows how to construct the twistor space $Z$ of a
GHC-manifold. Unfortunately, it also shows that $Z$ usually will not be
Hausdorff. Therefore we need, in general, to view the twistor space $Z$ as a foliation, given by the distribution \eqref{Inv}.
We adopt the following definition:
\begin{definition} The {\em twistor foliation} $Z$ of a generalised hypercomplex manifold $M$ is the foliation of $M\times {\Bbb C}P^1$ determined by the integrable distribution ${\cal Z}$, where ${\cal Z}_{(m,q)}= ({\cal F}_q)_m$. 
\end{definition}
\begin{definition} A GHC-manifold is {\em regular} if its twistor foliation is simple.\end{definition}
In the case of a regular GHC-manifold, we have a genuine twistor space which is a complex manifold fibering over $\cx P^1$. In general,
the twistor foliation $Z$ is a transversely holomorphic foliation. In other words, if we consider a small open subset of $U$ (where the foliation is simple), then the space of leaves of $U\times \cx P^1$ is a complex manifold. Moreover it fibers over $\cx P^1$. Globally, this means (directly from the definition) that the distribution ${\cal Z}$ is contained in the kernel of $d\pi$, where $\pi: M\times {\Bbb C}P^1\rightarrow {\Bbb C}P^1$ is the projection. We shall often abuse notation and write $\pi:Z\rightarrow {\Bbb C}P^1$.
We shall also refer to a leaf of $Z$ as a point of $Z$, and write $z\in Z$ (this is justified if we identify the space of leaves with a non-Hausdorff manifold). 
\par
We shall try as much as possible to limit ourselves to regular GHC-manifolds and just indicate how the results can be generalised.

\subsection{Complex GHC-manifolds}

If $M$ is a GHC-manifold arising as the space of real sections of a complex manifold $Z\rightarrow \cx P^1$, then $M$ has a natural complexification $M^\cx$ given as the space of all sections (not just real) with the correct normal bundle. For a general GHC-manifold we still have complexifications but not a canonical one.
We adopt the following definition:
\begin{definition} A {\em complex} GHC-manifold is a complex manifold $X$ with a holomorphic action of $SL(2,\cx)$ on $T^{1,0}X$ such that each tangent space is isomorphic to $S^k\cx^2\otimes \cx^n$, and such that the differential ideal generated by highest weight forms is closed for each Borel subgroup of  $SL(2,\cx)$. \end{definition}

As in the case of real GHC-manifolds we obtain, for each $q\in \cx P^1$,  a holomorphic subbundle ${\cal K}_q$ of $TX=T^{1,0}X$ and, putting these together, a holomorphic foliation $Z$ of $X\times \cx P^1$. Again, we call $Z$ the twistor foliation. Following the tradition, we shall call its leaves {\em $\alpha$-surfaces}.
\par
As before, we have the isomorphism:
$$ TX\simeq E_X\otimes H.$$
Consider the projection $\tau:X\times  \cx P^1\rightarrow X$ and the bundle $\tau^\ast E_X$. A simple open subset $U$ of $X\times  \cx P^1$ has a quotient manifold $U/{\cal K}$ which is a local twistor space. The proof of Proposition \ref{twistor} shows that $\tau^\ast E_X$ restricted to $U$ arises as a pullback of a bundle on  $U/{\cal K}$. Therefore:
\begin{proposition} The bundle $\tau^\ast E_X$ is a foliated bundle on $(X\times \cx P^1 ,Z)$.\hfill $\Box$\label{E_X}\end{proposition}

We observe that if $X$ is equipped with an antiholomorphic involution $\sigma$ compatible with the $SL(2,\cx)$-action, then the fixed point set $X^\sigma$ is a GHC-manifold, providing that $\dim_{\Bbb R}X^\sigma=\dim_\cx X$. In this case, the twistor foliations on $X$ and on $M=X^\sigma$ are {\em transversely equivalent} (\cite{Mol}, section 2.7) and so we can think of them as the same "manifold of leaves".

It will often be convenient to replace a GHC-manifold $M$ with a complex GHC-manifold $X$ such that $M=X^\sigma$. We shall then write $M^\cx$ for $X$.  


\section{On Quillen's resolution\label{Quillen}}

D. Quillen \cite{Q} has defined a canonical resolution of sheaves on ${\Bbb
C}P^1$. In this section we shall view this resolution and its splitting in terms
of homogeneous vector bundles on $\cx P^1$.
\par
Let $H$ be an irreducible representation of $SL(2,\cx)$ arising as the space of
sections of an ample line bundle $L={\cal O}(k)$ on ${\Bbb C}P^1$. We shall
write in this section $G$ for $SL(2,\cx)$ and $B$ for the standard Borel
subgroup of upper-triangular matrices. Thus $L$ is the homogeneous line bundle
$G\times_B \cx_k$, where $\cx_k$ is the $1$-dimensional representation of $B$:
$[b_{ij}]\cdot z=b_{11}^{-k}z$.
\par
We consider the homogeneous vector bundle $G\times_B H$.  Since the action of
$B$ on $H$ extends to an action of $G$, $G\times_B H$ is trivial as a vector
bundle. We have the canonical equivariant map
$$G\times_B H\rightarrow L$$
induced by the $B$-equivariant map $H\rightarrow \cx_k$. If we identify
$G\times_B H$ with the trivial bundle $\underline{ H}= H\times {\Bbb C}P^1$,
then this map simply sends $(h,q)$, $h\in H= \Gamma(L)$, to $h(q)$.
We obtain an exact sequence of homogeneous vector bundles on ${\Bbb C}P^1$:
\begin{equation}
0\rightarrow K \rightarrow G\times_B H \rightarrow L\rightarrow 0.
\label{K}\end{equation}
The cohomology sequence  of the dual to
\eqref{K} implies, given that $H^0(L^\ast)=H^1(\underline{H}^\ast)=0$,
\begin{equation}
0\rightarrow H^0(\underline{H}^\ast)\rightarrow H^0(K^\ast)\rightarrow
H^1(L^\ast)\rightarrow 0.\label{H0}
\end{equation}
This is an exact sequence of representations of $SL(2,\cx)$. The representation
$H^\prime=H^1({\Bbb C}P^1,L^\ast)$ is isomorphic to $S^{k-2}\cx^2$ and so
irreducible. Therefore   $\hat{H}=H^0({\Bbb C}P^1,K^\ast)$ is simply the direct
sum of $H^\ast$ and $H^\prime$ and the following sequence of representations
splits:
\begin{equation}
0\rightarrow H^\ast\stackrel{i}{\rightarrow} \hat{H}\stackrel{j}\rightarrow
H^\prime
\rightarrow 0.\label{split0}
\end{equation}
We observe:
\begin{lemma} The vector bundle $K^\ast$ splits as a direct sum of line bundles
${\cal O}(1)$.\label{O1}\end{lemma}
\begin{pf}
The long exact sequence of \eqref{K} implies that all cohomology of $K$ vanishes
(given that the second map induces an isomorphism on $H^0$), and hence $K$ is a
sum of ${\cal O}(-1)$'s.\end{pf}

Therefore $\hat H$ can be given structure of a complex quaternionic vector
space.
\par
We shall now consider the sequence \eqref{split0} in a greater detail and
construct a canonical equivariant isomorphism of  $H^1({\Bbb C}P^1,L^\ast)$ with
the space of sections of another homogeneous bundle.
 Let, for every $q\in {\Bbb C}P^1$, $S_q$ denote the subspace of highest weight
vectors in $H$ (for the Borel corresponding to $q$).
\begin{lemma} The bundle $S$ is a homogeneous subbundle of
$K$.\label{sub}\end{lemma}
\begin{pf} It is clear that $S$ is a homogeneous subbundle of $G\times_B H$,
since the map of fibers $S_1\rightarrow H$ over $[1]$ is $B$-equivariant. On the
other hand we also have $S_1\subset K_1$, since  $K_1=\{(v_1,\dots,v_k,0)^T\}$
and $S_1=\{(v_1,0,\dots,0)^T\}$.\end{pf}
We consider  the short exact sequence of homogeneous vector bundles:
\begin{equation} 0\rightarrow
(K/S)^\ast\rightarrow K^\ast\rightarrow S^\ast \rightarrow
0.\label{split2}\end{equation}The long
exact sequence of cohomology of this sequence begins as: $$0\rightarrow
H^0\bigl((K/S)^\ast\bigr)\rightarrow \hat{H}\rightarrow H^0(S^\ast).$$ The
construction of a line bundle from representation, recalled after Definition
\ref{GHC}, shows that there
is a canonical isomorphism (of representations) $H^\ast\rightarrow
H^0(S^\ast)$. Thus we obtain a map $p: \hat{H}\rightarrow H^\ast$.
\begin{lemma}
The map $p$ is a left inverse of the map $i$ in
\eqref{split0}.\label{split1}\end{lemma}
\begin{pf} This is a matter of going through various identifications. We start
with a section $s$ of the trivial bundle $\underline{H}^\ast$. The map $i$ means
that we evaluate $s$ on $K$. Following this by $p$ means that we evaluate $s$ on
$S$. This however is exactly how we obtain the element $s$ of $H^\ast$: it is a
section of $H^0(S^\ast)$.\end{pf}

Finally:
\begin{lemma} As homogeneous vector bundles
$$ (K/S)^\ast\simeq G\times_B H^\prime.$$
In particular, $K/S$ is a trivial vector bundle.\label{K/S}\end{lemma}
\begin{pf} As both sides are  homogeneous vector bundles, it is enough to show
that their fibers over $[1]$ coincide as $B$-modules. $H\simeq {\Bbb C}^{k+1}$
was the $k$-th symmetric power $S^k\cx^2$ of the standard representation of
$SL(2,\cx)$. The fiber $K_1$ of $K$ consists of vectors in $H$ of the form
$(v_1,\dots,v_k,0)^T$ and the fiber of $S_1$ of vectors $(v_1,0,\dots,0)^T$. Now
the mapping $K_1/S_1\rightarrow {\Bbb C}^{k-1}$ induced by:
$$(v_1,v_2\dots,v_k,0)^T\mapsto (v_2,\dots,v_k)^T$$
provides an isomorphism between $B$-modules $K_1/S_1$ and $S^{k-2}\cx^2$.
\end{pf}

\section{Some natural bundles on generalised hypercomplex
manifolds\label{bundles}}

Let $M$ be a regular GHC-manifold and $Z\rightarrow {\Bbb C}P^1$ its twistor
space. We have the double fibration
\begin{equation} Z\stackrel{\eta}{\longleftarrow} Y=M^\cx\times {\Bbb C}P^1
\stackrel{\tau}{\longrightarrow} M^\cx.\label{fibration}\end{equation}
The kernel of $d\eta$ is the bundle ${\cal K}$ (defined in section
\ref{definition}).
\par
 We shall consider several natural (locally free) sheaves over
$M^\cx$, which arise as direct images  of sheaves
on $Y$.
\par
We shall write $\tau^i_\ast(F)$ for the $i$-th direct image of the sheaf of
sections of $F$, and $\tau_\ast$ for $\tau_\ast^0$. The first bundle on
$M^\cx$ we wish to consider is the tangent bundle $TM^\cx=T^{1,0}M^\cx$. This,
by
the proof of Proposition \ref{twistor} can be written as $\tau_\ast\eta^\ast
(T_\pi Z)$, where $T_\pi Z$ is the vertical tangent bundle. We recall the decomposition \eqref{E}: \begin{equation} TM^\cx\simeq E_M\otimes H.\label{tensor}\end{equation}
The proof of  Proposition \ref{twistor} shows that $E_M= \tau_\ast\eta^\ast\bigl(T_\pi Z\otimes \pi^\ast L^\ast\bigr)$
and $H$ is the trivial bundle with fiber equal to the representation
$H^0({\Bbb C}P^1,L)=S^k\cx^2$ of $SL(2,\cx)$. If $k$ is odd, so that $L$ has a
quaternionic structure, then both $E_M$ and $H$ have quaternionic structures
covering the real structure on $M^\cx$. Similarly, if $k$ is even, then both
$E_M$ and $H$ have real structures. In particular, if $k$ is even, the bundle
$E$ has a real form $E^{\Bbb R}$ on $M$, and $TM\simeq E^{\Bbb R}\otimes H^{\Bbb
R}$, where $ H^{\Bbb R}$ is the underlying representation of $SU(2)$ on ${\Bbb
R}^{k+1}$.
\par
Next we consider direct images of $\eta^\ast (T^\ast_\pi Z)$.
\begin{proposition} We have: $$\tau^0_\ast \eta^\ast\left( T^\ast_\pi Z\right)=
0$$ $$\tau^1_\ast \eta^\ast\left( T^\ast_\pi Z\right)\simeq E^\ast \otimes
H^\prime, $$ where $E=E_M$ and $H^\prime$ is the trivial bundle over $M^\cx$
with fiber given by the representation $S^{k-2}\cx^2\simeq H^1({\Bbb
C}P^1,L^\ast)$ of $SL(2,{\Bbb C})$.\label{dir}
\end{proposition}
\begin{pf} We can write $T^\ast_\pi Z\simeq \bigl(T_\pi^\ast Z\otimes
\pi^\ast(L)\bigr)\otimes \pi^\ast (L^\ast)$. Since the first bundle in the
tensor product is trivial on each twistor line and $L$ is positive, the result
follows.\end{pf}

Now we consider the sheaf $\Omega_\eta^\ast$ of $\eta$-vertical
holomorphic forms  on $Y$, i.e. the exterior algebra of  
$\Omega^1(Y)/\eta^\ast\bigl(\Omega^1(Z)\bigr).$ We
are particularly interested in the direct image $\tau_\ast \Omega^1_\eta$
(it is well known \cite{W,Ma,B-E} that the Ward correspondence gives us
differential operators $F\rightarrow \tau_\ast\Omega^1_\eta\otimes F$ on
bundles over $M$).

We have
\begin{proposition}
The sheaf $\tau_\ast (\Omega^1_\eta)$ is isomorphic to $E^\ast \otimes
\hat{H}$, where $\hat{H}=H^\ast\oplus H^\prime$ is the representation defined in
the previous section.\label{phi}
\end{proposition}
\begin{pf}
Let $x\in M^\cx$ and let ${\Bbb C}P^1_x$ be the fibre of $\tau$ over $x$ (i.e. a
section of $Z$ corresponding to $x$). The $\eta$-normal bundle of ${\Bbb
C}P^1_x$
inside $Y$ is the bundle ${\cal K}$, i.e. the subbundle of $T_x\times {\Bbb
C}P^1$
given by choosing the subspace ${\cal K}_q$ of $T_x$ for each $q$. Therefore
\begin{equation}\left(\tau_\ast \Omega^1_\eta\right)_x= \Gamma({\Bbb
C}P^1_x,{{\cal K}}^\ast).\label{omega}\end{equation}
 However, as $T_x$ decomposes as ${\Bbb
C}^n\otimes H$, where ${\Bbb C}^n$ is the fiber of $E$ at $x$,  the bundle
${\cal K}$ decomposes as ${\Bbb C}^n \otimes K$, where $K$ is defined in
\eqref{K}. The result follows.
\end{pf}

We observe that $\Omega^1_\eta$ fits into the exact sequence
\begin{equation} 0\rightarrow  \eta^\ast (T^\ast_\pi Z )\rightarrow \tau^\ast
\left(\Omega^1M^\cx\right) \rightarrow \Omega^1_\eta\rightarrow 0.
\label{eta}\end{equation} The corresponding long exact sequence of direct
images gives,  using Proposition \ref{dir}:
\begin{equation} 0\rightarrow \Omega^1M^\cx \rightarrow \tau_\ast\Omega^1_\eta
\rightarrow E^\ast\otimes H^\prime\rightarrow 0.\label{split}\end{equation}

All the sheaves in this sequence are tensor products of $E^\ast$ with
(sheaves of sections of) trivial bundles and we have the key:

\begin{proposition} On a regular GHC-manifold the sequence \eqref{split} splits
canonically:
\begin{equation}\tau_\ast \Omega^1_\eta=\Omega^1\oplus \left( E^\ast\otimes
H^\prime\right).\end{equation}
 The splitting is induced by that of Lemma \ref{split1}.
\label{splits}\end{proposition}
\begin{pf}
The only thing to prove is that the maps in \eqref{split} are identity on
$E^\ast$. This follows from the identification of the bundles involved at
each point $x$ of $M^\cx$, as in the proof of Proposition \ref{phi}.
\end{pf}

It is now clear how to compute direct images of $\Omega_\eta^i$ for any $i$. For
example,  the same arguments as in the proof of Proposition \ref{phi} yield:
\begin{proposition}
The sheaf $\tau_\ast (\Omega^2_\eta)$ is isomorphic to $\bigl(S^2E^\ast \otimes
H_-\bigr)\oplus\bigl(\Lambda^2E^\ast\otimes H_+\bigr)$, where $H_-=H^0({\Bbb
C}P^1,\Lambda^2 K^\ast)$ and $H_+=H^0({\Bbb C}P^1,S^2 K^\ast)$.\hfill
$\Box$.\label{omega2}\end{proposition}

\begin{remark} Although we have limited ourselves to regular GHC-manifolds, everything in this section remains valid for arbitrary GHC-manifolds, as the definitions of sheaves $\Omega^\ast_\eta$ and of $\eta^\ast(T^\ast_\pi Z)$ ("sheaf of vertical basic $1$-forms") are unchanged. 
\end{remark}

\section{Monopoles\label{monopoles}}

We shall consider geometric structures arising from holomorphic vector
bundles $F$  on the twistor space $Z$ of a regular generalised hypercomplex
manifold $M$, under
the condition that $F$ is trivial on twistor sections.
We use the notation of the previous section, in particular the double fibration
 $$ Z\stackrel{\eta}{\longleftarrow} Y=M^\cx\times \cx P^1
\stackrel{\tau}{\longrightarrow} M^\cx,$$
and the sheaf  $\Omega_\eta^\ast$ of $\eta$-vertical holomorphic forms. By composing the exterior derivative on $Y$ with the projection onto $\Omega_\eta^1 $ we obtain a first-order differential operator 
$$ d_\eta: \Omega^0\rightarrow \Omega_\eta^1$$
which annuls $\eta^\ast  \Omega^0_Z$. If $F$ is  a holomorphic bundle
on $Z$ then it is well known (see, e.g. \cite{Ma} or \cite{B-E}) that $d_\eta$  extends to a flat relative connection
$$\nabla_\eta: \eta^\ast F\rightarrow \Omega_\eta^1(F)$$
on the pullback $\eta^\ast F$. 
\par
 If the bundle $F$ is trivial on each section of $Z$, then
$\eta^\ast F$ is trivial on each fibre of $\tau$ and we consider the direct
image $$\hat{F}=\tau_\ast\eta^\ast F$$ which is a vector bundle on $M^\cx$
of the same rank as $F$. By pushing down $\nabla_\eta$ we obtain a
first-order differential operator
\begin{equation} D:=\tau_\ast\nabla_\eta:\hat{F}\rightarrow \tau_\ast
\Omega^1_\eta\otimes\hat{F}.\end{equation}
This operator satisfies:
$$D(fs)=fD(s)+\partial f\otimes s,$$
where $\partial=\tau_\ast d_\eta$.
\par
We recall from the previous section that we have a canonical isomorphism
\begin{equation} \left(\tau_\ast
\Omega^1_\eta\right)_x =H^0(\cx P^1_x,{\cal K}^\ast),\end{equation}
where the bundle ${\cal K}$ is the subbundle of $T_x\times \cx P^1$ annihilated by the
highest weight forms for each $q\in \cx P^1$. Therefore we have the canonical
map
\begin{equation} \imath_q^\ast:\tau_\ast \Omega^1_\eta\rightarrow
{{\cal K}}_q^\ast\label{P}\end{equation} for each $q\in \cx P^1$. In
particular, if we restrict both $\hat{F}$ and $D$ to the leaf of the
foliation ${\cal K}_q$ passing through $x$, i.e. to the submanifold
$\tau\bigl(\eta^{-1}(z)\bigr)$ where $z=\eta(\tau^{-1})(x)$, then we recover
$\nabla_\eta$, i.e. $D$ induces a flat connection on $\hat{F}$ restricted to
this submanifold.
\par
All of the above is well-known and works for arbitrary twistor spaces
\cite{Ma,B-E}. In the case of generalised hypercomplex manifolds we can use,
however, the
results of the previous section, in particular the splitting:
\begin{equation}\tau_\ast \Omega^1_\eta=\Omega^1\oplus \left( E^\ast\otimes
H^\prime\right).\label{reallysplits}\end{equation}
The operator $\partial$ becomes simply the exterior derivative:
\begin{lemma} Under the isomorphism \eqref{reallysplits}, the operator
$\partial: \Omega^0\rightarrow \tau_\ast \Omega^1_\eta$ becomes
$\partial=d\oplus 0$.\end{lemma}
\begin{pf}
From the equation \eqref{eta} we have the commutative diagram:
$$\begin{CD} 0 @>>> \Omega^0 @= \Omega^0 @>>> 0\\
 @. @V{d}VV  @V{\partial}VV @. \\
0 @>>> \Omega^1 @>>> \tau_\ast\Omega^1_\eta
@>>> E^\ast\otimes H^\prime.\end{CD}$$
\end{pf}

Therefore on a generalised hypercomplex manifold the operator $D$ is given
by a connection $\nabla$  and a Higgs field, i.e. a section $\Phi$ of
$\End(\hat{F})\otimes (E^\ast \otimes H^\prime)$.
\par
Moreover $\nabla+\Phi$ is flat on each $\alpha$-subspace, i.e. on each
submanifold $S_z$ of  $M^\cx $ of the form $\tau\bigl(\eta^{-1}(z)\bigr)$,
$z\in Z$.

We can formulate the results as follows:
\begin{theorem} Let $M$ be a regular generalised hypercomplex manifold and $Z$
its twistor space. There exists a 1-1 correspondence between
\begin{itemize}
\item[(a)] holomorphic vector bundles on $Z$
trivial on each  twistor line, and
\item[(b)] {\em monopoles} on $M^{\Bbb C}$, i.e. a connection $\nabla$ on a
holomorphic vector bundle $\hat{F}$ on $M^\cx $ and a section $\Phi$ of
$\End(\hat{F})\otimes \left(E^\ast\otimes H^\prime\right)$ such that
$\nabla\oplus \Phi$ is flat on each $\alpha$-subspace, i.e. on each
submanifold of  $M^\cx $ of the form $\tau\bigl(\eta^{-1}(z)\bigr)$, $z\in
Z$.\end{itemize} This correspondence remains valid in the presence of a real
structure, giving monopoles on $M$.\hfill $\Box$\label{MONO}\end{theorem}
Here, the connection $\nabla\oplus \Phi$ on an $\alpha$-subspace $S_z$
corresponding to a point $z\in Z$ over $q\in \cx P^1$ is given by: $$
\hat{F}\stackrel{\nabla\oplus \Phi}{\longrightarrow} \left(\hat{F}\otimes
E^\ast\otimes H^\ast\right)\oplus \left(\hat{F}\otimes E^\ast\otimes
H^\prime\right) =\hat{F}\otimes E^\ast\otimes \hat{H}
\stackrel{\imath^\ast_q}{\longrightarrow} \hat{F}\otimes \Omega^1S_z,$$
where $\imath^\ast_q$ is given by \eqref{P}.

\begin{remark}
 The above theorem describes monopoles for the group $GL(m,\cx)$ (or  $GL(m,{\Bbb
R})$), $m$ being the rank of $F$. By considering bundles whose structure
group reduces we obtain a 1-1 correspondence between bundles on $Z$ and
monopoles for other groups $G$. \end{remark}

\begin{example}({\bf Monopoles on ${\bf\Bbb R}^5$})  Let $M$ be ${\Bbb R}^5$
viewed as the real form of the fourth symmetric power of the defining representation of ${\frak
su}(2)$. Then $M$ is a $4$-hypercomplex manifold and its twistor space is the
total space of the line bundle ${\cal O}(4)$. Thus $M^\cx$ is identified with
polynomials
\begin{equation} z_0+z_1\zeta+z_2\zeta^2+z_3\zeta^3 +z_4\zeta^4,
\label{sect}\end{equation}
and $M$ with polynomials invariant under \begin{equation} z_i\mapsto
(-1)^i\overline{z}_{4-i}.\label{reality}\end{equation}
\par
We wish to discuss the monopoles on a (trivial) vector bundle $F$ over ${\Bbb
R}^5$ (cf. \cite{MS}). We first discuss monopoles on ${\Bbb C}^5$. In this case
$E_M$ is the trivial $1$-dimensional bundle and, according to the above theorem,
a monopole is given by a connection and and a section $\Phi$ of $\End(F)\otimes
H^\prime\simeq \End(F)\otimes {\Bbb C}^3$, i.e. a triple $\Phi_1,\Phi_2,\Phi_3$
of sections of $\End(F)$. Let the $1$-form of the connection be
$$\sum_{i=0}^4 A_idz_i.$$
In order to find the monopole equation we have to consider the map
$\imath^\ast_q:{\Bbb C}^5\oplus {\Bbb C}^3\rightarrow K_q^\ast$ given by
\eqref{P}. Its image consists of polynomials \eqref{sect} vanishing at
$\zeta=q$, and $\imath^\ast_q$ on ${\Bbb C}^5$ is simply the projection. On the
other hand $\imath^\ast_q$ on ${\Bbb C}^3$ is described by Lemma \ref{K/S} and the equation \eqref{split2}.
It is equivariant for the Borel subgroup corresponding to $q$.
\par
Consider $q=0$. Then the connection $\nabla^0$ on each $\alpha$-surface
$S_z$ obtained by freezing $z_0$ in \eqref{sect} has the $1$-form: $$
(A_1+\Phi_1)dz_1+(A_2+\Phi_2)dz_2+(A_3+\Phi_3)dz_3+A_4dz_4.$$ Thus this
connection must be flat and similarly for every $q$.  We would obtain $18$
equations due to Mason and Sparling \cite{MS}. In practice it is simpler to
obtain the equations as a reduction of self-duality equations on the
connection $$\sum_{i=1}^3\Phi_i dt_i+\sum_{i=0}^4 A_idz_i$$ on ${\Bbb R}^8$.
This will be discussed in detail in section \ref{ext}.
\par
 Monopoles on ${\Bbb R}^5$
will arise when we impose the reality condition \eqref{reality} on the $A_i$
and a similar one on the $\Phi_i$.
\end{example}

We briefly discuss the (well-known) case of $1$-hypercomplex manifolds. In
this case there is no $\Phi$ and a connection $\nabla$ obtained as in
Theorem \ref{MONO} is called {\em self-dual} \cite{CS} or {\em hyperholomorphic} \cite{V}. It is easy to
describe this condition. Since, for $k=1$, $\tau_\ast \Omega^1_\eta=
\Omega^1$, we have the natural operator
\begin{equation}
 \Omega^2=\tau_\ast\Omega^1_\eta\wedge \tau_\ast\Omega^1_\eta
\longrightarrow \tau_\ast\Omega^2_\eta.\label{curv}\end{equation}
The curvature of the relative connection $\nabla_\eta$ lies in
$\Omega^2_\eta$, and it vanishes. Therefore a connection $\nabla$ will be
self-dual if and only if its curvature lies in the kernel of \eqref{curv}
\cite{B-E}. Proposition \ref{omega2} describes $\tau_\ast\Omega^2_\eta$ and,
since in this case $K^\ast={\cal O}(1)$ and $H=\cx^2$, we conclude:
\begin{proposition}
A connection $\nabla$ on a holomorphic vector bundle over the complexified
hypercomplex manifold $M^\cx$ is self-dual if and only if its curvature lies
in the component $S^2E^\ast \otimes \Lambda^2 H$ of $\Lambda^2T^\ast M\simeq
\bigl(S^2E^\ast \otimes \Lambda^2 H\bigr)\oplus \bigl(\Lambda^2E^\ast
\otimes S^2 H\bigr)$.\hfill $\Box$
\end{proposition}

An equivalent way of expressing this condition is that the curvature is
$SL(2,{\Bbb C})$-invariant \cite{V}.

\begin{remark} Again, the results of this section hold for non-regular GHC-manifolds, providing we replace holomorphic bundles on $Z$ with foliated holomorphic bundles on $(Y,Z)$ (cf. definition \ref{f-bundle}). In particular, there is a monopole on the bundle $E_M$ over $M^\cx$ .\end{remark}

\begin{remark} Let $M$ be a $k$-hypercomplex manifold of dimension $n(k+1)$ ($n$ is even if $k$ is odd). Since the bundle $E$ has a natural connection $\nabla^E$ and $T^\cx M=E\otimes H$, we obtain a canonical connection $\nabla_M$ on $M$ by tensoring $\nabla^E$ with the flat connection on $H$. Since the action of $SU(2)$ on $TM$ is parallel for this connection, the holonomy of $\nabla_M$ is contained in the centraliser of $SU(2)$ in $GL(n(k+1),{\Bbb R})$, i.e. in $GL(n,{\Bbb R})$, if $k$ is even and in $GL(n/2,{\Bbb H})$ if $k$ is odd. Apart from $k=1$, the connection $\nabla_M$ is not torsion-free. It is perhaps of interest to know more about $\nabla_M$, e.g. how its torsion is related to the Higgs field arising from $E$.\end{remark}

\section{Structure of generalised hypercomplex manifolds\label{str}}

From the very definition,  a hypercomplex manifold has the structure of a
complex manifold for every $q\in {\Bbb C}P^1$. For generalised hypercomplex
manifolds, we have shown so far (in section \ref{definition}) that a 
GHC-manifold $M$ has,  for every $q\in {\Bbb C}P^1$, the structure of a
transverse holomorphic foliation $(M,Z_q)$, where leaves of $Z_q$ are the integral manifolds of the distribution ${\cal F}_q$. In this section we shall prove a much more precise structure theorem.

\subsection{ Weight subbundles}
Let $M$ be an almost $k$-hypercomplex manifold. Then, for every $q\in \cx P^1$,
we can decompose $T^\cx M$ into eigenspaces of the circle normalising the Borel
subgroup corresponding to $q$:
\begin{equation} T^\cx M=\bigoplus_{i=0}^{k} {\cal
S}_q(k-2i),\label{decomposition}\end{equation}
where we adopt the convention that ${\cal S}_q(-k)$ is the bundle of lowest
weight vectors and ${\cal S}_q(i+1)$ is obtained from ${\cal S}_q(i)$ by the
action of the unipotent radical of $B_q$. In particular the subbundle ${\cal
K}_q$ defined in section \ref{definition} is the bundle $\bigoplus_{j\neq -k}
{\cal S}_q(j)$.
\par
On a $k$-hypercomplex manifold these distributions are defined on all of
$M^\cx$. We have:
\begin{theorem}
The subbundles   ${\cal S}_q(j)$ of $TM^\cx$ are involutive for every $j$ and
$q$.\label{S^j}\end{theorem}

This follows from the following result:
\begin{proposition} For every $q\in \cx P^1$ and every $l\geq -k$, the subbundle
$$ {\cal K}_q(l)=\bigoplus_{j\geq l}  {\cal S}_q(j)$$
of $TM^\cx$ is involutive.\label{K(l)} \end{proposition}
We observe that ${\cal K}_q(-k)=TM^\cx$ and ${\cal K}_q(-k+1)={\cal K}_q$, where ${\cal K}_q$ is the leaf of the twistor foliation over $q$. Theorem \ref{S^j} follows immediately from this proposition, as ${\cal S}_q(j)=  {\cal K}_q(j)\cap  {\cal K}_{\sigma(q)}(-j)$, where $\sigma$ is the antipodal map on $\cx P^1$.
\par
To prove the above proposition, consider a regular neighbourhood $U$ of a point $m\in M^\cx$ with the corresponding twistor space $Z_U$. Then the integral submanifold of $ {\cal K}_q(l)$ through $m$ is the space of all sections of $Z_U\rightarrow \cx P^1$ which coincide with the section $m:\cx P^1\rightarrow Z_U$ up to order $l-1$ at $q$ (equivalently these are sections of $Z_U$ blown up $l-1$ times lying in the same connected component as $m$).

\subsection{Structure theorems}

We shall now show that, for every $q\in \cx P^1$,  a generalised hypercomplex manifold $M$ {\em locally} looks  like the following sequence of submersions:
$$M \rightarrow Z_q(p)\rightarrow \dots \rightarrow Z_q(1)=Z_q,$$
where $p=(k+1)/2$ if $k$ is odd and $p=k/2$ if $k$ is even. Each $Z_q(i)$ is a complex manifold of dimension $2ni$ if $k$ is odd and $ni$ if $k$ is even. Moreover the projections $Z_q(i)\rightarrow Z_q(i-1)$ are holomorphic. 
\par
We first define the relevant pseudogroups of transition functions. For $k$ odd, let $\Gamma_{k,n}$ be the pseudogroup of local holomorphic diffeomorphisms of ${\Bbb C}^{2n}\otimes\cx^{(k+1)/2}$ of the form 
$$ (x_1,x_2,\dots,x_{\frac{k+1}{2}})\longmapsto \bigl(\phi_1(x_1),\phi_2(x_1,x_2), \phi_3(x_1,x_2,x_3),\dots, \phi_{\frac{k+1}{2}}(x_1,\dots,x_{\frac{k+1}{2}})\bigr),$$
where for every  $i>0$, $\phi_i$ is a local holomorphic transformation of ${\Bbb C}^{2n}$. \newline
For $k$ even, we define $\Gamma_{k,n}$ similarly as the pseudogroup of local diffeomorphisms of $\bigl({\Bbb C}^{n}\otimes\cx^{k/2}\bigr)\oplus {\Bbb R}^n$, where each $\phi_i$ for $i=1,\dots,\frac{k}{2}$ is as before and $\phi_{\frac{k}{2}+1}(x_1,\dots,x_{\frac{k}{2}}, x_{\frac{k}{2}+1})$ is real in the "last" variable.

We now have the following structure theorems:

\begin{theorem}
Let $k$ be odd and let $M$ be a $k$-hypercomplex manifold of dimension
$2n(k+1)$. Then, for every $q\in \cx P^1$, $M$ has a complex atlas with transition functions belonging to $\Gamma_{k,n}$.\label{str1}\end{theorem}

\begin{theorem} Let $k$ be even and let $M$ be $k$-hypercomplex
manifold of dimension $n(k+1)$. Then, for every $q\in \cx P^1$,  $M$ has an atlas of charts in $\cx^{nk}\oplus {\Bbb R}^n$  with transition functions belonging to $\Gamma_{k,n}$.
\label{str2}\end{theorem}

\begin{pf} These are a simple consequence of Proposition \ref{K(l)} and the Newlander-Nirenberg theorem \ref{NN}.\end{pf}

\begin{remark} For $k$ odd, $M$ becomes a complex manifold. The proof shows how
to describe the complex structure: on each weight subbundle $F_i=\bigl({\cal
S}_q(i)\oplus {\cal S}_q(-i)\bigr)\cap TM$ of $TM$ choose an element $I_i$ of
the circle normalising $B_q$ such that  $I_i^2=-1$ on $F_i$. Then $I_q=\bigoplus
I_i|F_i$ is a complex structure on $M$.\end{remark}

\section{From $k$-hypercomplex to $(k-2i)$-hypercomplex\label{blow}}

We consider again the decomposition \eqref{decomposition} of $T^\cx M$ into weight subbundles. For $q\in {\Bbb C}P^1$ and any $l>0$, $l<k/2$, we consider the subbundle of $TM$ given by:
\begin{equation} {\cal F}_q(l)=TM\cap \bigoplus_{i=0}^{k-2l}  {\cal S}_q(k-2l-2i).\end{equation}
Thus ${\cal F}_q(0)=TM$ and ${\cal F}_q(1)={\cal F}_q$, where ${\cal F}_q$ is the leaf of the twistor foliation over $q$. 
\par
We are going to sketch a proof (cf. \cite{BG, DM, DM2}) of:
\begin{theorem} Let $M$ be a $k$-hypercomplex manifold. Then, for any $q\in {\Bbb C}P^1$ and any $0\leq l\leq k/2$,  every leaf $L$ of $ {\cal F}_q(l)$ carries a canonical $(k-2l)$-hypercomplex structure. In particular, if $k$ is odd, a  $2n(k+1)$-dimensional  $k$-hypercomplex manifold $M$ is foliated by $4n$-dimensional hypercomplex manifolds.\label{k-2}\end{theorem}
\begin{pf} We first construct an almost $(k-2l)$-structure on each leaf $L$. 
The complexified tangent space at $m$ to $L$ is canonically the subspace of $T_m^\cx M$ consisting of all but $l$ highest and $l$ lowest weight subspaces for the Borel $B_q$. This subspace has canonically the structure of  an $SL(2,\cx)$-module with highest weight $k-2l$ (just set the actions of $B_q$ and of the opposite Borel to be zero on $S_q(k-2l)$ and on $S_q(-k+2l)$, respectively). As everything is compatible with the real structure, we obtain an almost $(k-2l)$-hypercomplex structure on $L$. To prove its integrability, consider a small neighbourhood $U$ of $m$, where the twistor foliation is simple. Let $Z_U$ be the twistor space of $U$. Consider first the case $l=1$. A leaf of  ${\cal F}_q(1)$ is simply a fiber $\eta^{-1}(z)$ of $\eta: U\rightarrow Z_U$, where $q=\pi(z)$. Its complexification is the intersection of the fibers $\eta^{-1}(z)$ and $\eta^{-1}(\tau(z))$ of $\eta: U^\cx\rightarrow Z_U$, where $\tau:Z_U\rightarrow Z_U$ is the real structure. One can check that this intersection is the space of sections passing through the exceptional divisors of the manifold $\tilde{Z}_U$ which is $Z_U$ blown up at $z$ and $\tau(z)$. The almost $(k-2)$-hypercomplex structure on $L$ is the one obtained from Proposition \ref{twistor} applied to $\tilde{Z}_U$ and hence it is integrable. Blowing up successively proves the result.
\end{pf}

\section{A hypercomplex extension of a $k$-hypercomplex manifold\label{ext}}

If $M$ is a flat $k$-hypercomplex manifold, i.e. a vector space of the form
$H^{\Bbb R}\otimes{\Bbb R}^n$, where  $H^{\Bbb R}$ is the real form of
$S^k\cx^2$, then the sequence dual to \eqref{split0} can be interpreted as
saying that there is a $1$-hypercomplex manifold $\tilde{M}$ fibering over $M$
and such that $T\tilde{M}$ can be canonically identified with
$\tau_\ast\Omega^1_\eta$. In particular monopoles on $M$ will correspond to
self-dual connections on  $\tilde{M}$. In this section we shall show that such
an $\tilde{M}$ always exists.
Let $M$ be a $k$-hypercomplex manifold. Then $TM^\cx\simeq E_M\otimes H$ and
both $E_M$ and $H$ are equipped with real or with quaternionic structures. The
same is true then about $E_M$ and $H^\prime\simeq S^{k-2}\cx^2$ and we have a
real vector bundle $(E_M\otimes H^\prime)^{\Bbb R}$ over $M$.
\par
The total space of $(E_M\otimes H^\prime)^{\Bbb R}$ is an obvious candidate for $\tilde{M}$. One way of defining an almost hypercomplex structure is to use the connection on $E_M$ to define a connection on    $\tilde{M}^\cx=E_M\otimes H^\prime$ and obtain a splitting 
\begin{equation}T\tilde{M}^\cx = p^\ast (E_M\otimes H^\prime)\oplus  p^\ast (E_M\otimes H)\simeq  p^\ast (E_M\otimes \hat{H}),\label{twist}\end{equation}
where $p:\tilde{M}^\cx\rightarrow M^\cx$ is the projection.
Since $\hat{H}=H^0(K^\ast)$, Lemma \ref{O1} provides an almost hypercomplex structure on $\tilde{M}$. This will not, however, usually be integrable. We have to twist the first term of the middle expression in \eqref{twist} by the monopole on $E_M$ in order to get an integrable hypercomplex structure. Such constructions are given in \cite{BC}. The hypercomplex structure exists only on an open subset where the monopole is non-degenerate. 
We shall give here a twistorial and local proof that such an integrable structure exists, providing that certain topological conditions are fulfilled. 
To define these, consider the bundles ${\cal S}_q(k)$ of highest weight vectors on $M^\cx$. These
bundles combine to give us a bundle ${\cal S}$ on $Y=M^\cx\times \cx P^1$ which is a subbundle of
${\cal Z}$, i.e. of the twistor distribution. From the previous section we know that the bundles ${\cal S}_q(k)$ are involutive, and, hence, so  is  ${\cal S}$. \newline
We can state the theorem: 
\begin{theorem}
Let $M$ be a regular $k$-hypercomplex manifold such that the space of leaves of the distribution ${\cal S}$ is a manifold  $\tilde{Z}$.  Then there exists a hypercomplex manifold $\tilde{M}$ such that:
\begin{itemize}
\item There is a projection $p:\tilde{M}\rightarrow M$ and a section $M\rightarrow \tilde{M}$ of this projection, whose tubular neighbourhood can be identified with a  neighbourhood of the zero section of the bundle $(E_M\otimes H^\prime)^{\Bbb R}$ over $M$.
\item There is a canonical identification of $\Omega^1_{\tilde{\eta}}$ with $p^\ast\Omega^1_\eta$ of sheaves on $\tilde{Y}=\tilde{M}^\cx\times \cx P^1$, which makes the following diagram commute 
\begin{equation}
\begin{CD}\Omega^1_{\tilde{Y}} @>>> \Omega^1_{\tilde{\eta}}\\
 @AAA   @|  \\
p^\ast\Omega^1_Y @>>> p^\ast\Omega^1_\eta. \end{CD}\label{last}\end{equation}\end{itemize}
\label{EXT}\end{theorem}
\begin{corollary} In the above situation, there is a canonical isomorphism
\begin{equation} \Omega^1\tilde{M}^\cx\simeq p^\ast\tau_\ast \Omega^1_\eta. \label{identify}\end{equation}
\end{corollary}
\begin{pf}
Under the assumptions we have a well defined twistor space $Z$ of $M$ and the Hausdorff manifold $\tilde{Z}$ of leaves of ${\cal S}$. We shall show that  $\tilde{Z}$ is a twistor space of a hypercomplex manifold. 
Consider the double fibration  $ Z\stackrel{\eta}{\longleftarrow} Y=M^\cx\times
\cx P^1
\stackrel{\tau}{\longrightarrow} M^\cx$.  Since  ${\cal S} \subset   \Ker d\eta$, $\tilde{Z}$ fibers over $Z$  and has a canonical real structure given by
the real structures on $Y$ and on ${\cal S}$. A point of $M^\cx$ gives a section of $Y\rightarrow \cx P^1$ and hence a section of $\tilde{Z}$. We claim the normal bundle of such a section in  $\tilde{Z}$ splits into the direct sum of ${\cal O}(1)$'s. Indeed, the vertical tangent
bundle of the fibration $Y\rightarrow \cx P^1$ is simply  $\tau^\ast E_M\otimes
\underline{H}$, where $\underline{H}$ is the trivial bundle over $\cx P^1$ whose
fiber is $H=S^k\cx^2$. Therefore the vertical tangent bundle of $\tilde{Z}$ is
$$ T_{\tilde{\pi}}\tilde{Z}=\bigl(\tau^\ast E_M\otimes \underline{H}\bigr)/{\cal
S}\simeq \tau^\ast E_M\otimes (\underline{H}/S),$$
where $S$ is the bundle of highest weight vectors in $\underline{H}$.
The sequences \eqref{K} and \eqref{split2} imply that  $\underline{H}/ S\simeq
K^\ast$ and by Lemma \ref{O1} this splits into the direct sum of ${\cal O}(1)$'s. Since $\tau^\ast E_M$ is trivial on sections, we conclude that the space of real sections of $\tilde{Z}$ is a hypercomplex manifold, at least in a neighbourhood $\tilde{M}$ of $M$.
\par
As any section of $\tilde{Z}$ projects to a section of $Z$, we obtain a projection $p:\tilde{M}\rightarrow M$. To identify a tubular neighbourhood of $M$ in $\tilde{M}$ observe that   $T_{\tilde{\pi}}\tilde{Z}$ fits into
the diagram:
$$\begin{CD}0\longrightarrow  (T_\eta Y)/{\cal S} @>>>
T_{\tilde{\pi}}\tilde{Z} @>>>T_\pi Z @>>> 0\\
 @| @| @| @. \\
0\longrightarrow\tau^\ast E_M\otimes (K/S) @>>>   \tau^\ast E_M\otimes
(\underline{H}/S)@>>> (T_\pi Z\otimes L^\ast)\otimes L  @>>> 0.\end{CD}$$
This allows, using Lemma \ref{K/S}, to identify a tubular neighbourhood of $M^\cx$ in $\tilde{M}^\cx$ with a neighbourhood of the zero section in the bundle $E_M\otimes H^\prime$ over $M^\cx$.
\par
We now prove the second statement. From the above diagram, 
$T\tilde{Y}$ is identified with $(p^\ast\tau^\ast E_M)\otimes \hat{H}^\ast$ and 
$\tilde{\eta}^\ast \Omega^1_{\tilde{\pi}}$ with $(p^\ast\tau^\ast E^\ast_M)\otimes (H/S)^\ast$. Therefore  $\Omega^1_{\tilde{\eta}}$ is  $p^\ast \tau^\ast E^\ast_M$ tensored with the cokernel of the map $( H/S)^\ast\rightarrow\hat{H}$. Homogeneous bundle arguments, similar to those in section \ref{Quillen}, show that this cokernel is equal to $K^\ast$. Hence $\Omega^1_{\tilde{\eta}}$ is  isomorphic to $(p^\ast \tau^\ast E^\ast_M)\otimes K^\ast$ which is $p^\ast \Omega^1_\eta$ (see section \ref{bundles}). The diagram \eqref{last} commutes, as the  maps are identity on $p^\ast\tau^\ast E_M$ and the horizontal maps are given by appropriate maps on sheaves on $\cx P^1$, as in section \ref{Quillen}.
\end{pf}

The above theorem allows us to view monopoles on $M$ as a reduction of self-dual
connections (i.e. monopoles without $\Phi$) on $\tilde{M}$. Indeed, we recall from section \ref{monopoles} that a monopole is equivalent to a first order operator  $D:{F}\rightarrow \tau_\ast
\Omega^1_\eta\otimes{F}$ on a bundle $F$. Pulling $D$ back to $\tilde{M}^\cx$ we obtain a connection $\tilde{\nabla}$ on $p^\ast F$ (notice that \eqref{last} guarantees that the isomorphism   $\Omega^1_{\tilde{\eta}}\simeq p^\ast\Omega^1_\eta$ commutes with differentials $d_{\tilde{\eta}}$ and $d_\eta$). We have:

\begin{theorem} Let $M$ and 
$\tilde{M}$ be as in the previous theorem.
Then a pair  $(\nabla, \Phi)$ on a bundle $F$ over $M^\cx$ is a monopole if and only
if the connection  $\tilde{\nabla}$ on  $p^\ast F$ over  $\tilde{M}^\cx$ is
self-dual.\end{theorem}
\begin{pf} This is now automatic, since 
both $(\nabla, \Phi)$ and $\tilde{\nabla}$ arise from flat relative connections on $\Omega^1_\eta$ and $\Omega^1_{\tilde{\eta}}$, respectively. \end{pf}

\section{Maps between generalised hypercomplex manifolds\label{maps}}

We wish to consider maps between GHC-manifolds. What we clearly need are maps
which, for regular GHC-manifolds, give rise to fibrewise mappings of the twistor spaces. The following
definition is the translation of this condition.
\begin{definition}
 A morphism between two GHC-manifolds is a smooth map $f:M\rightarrow M^\prime$
such that, for any $q\in {\Bbb C}P^1$, the differential $df$
 satisfies the following two conditions
 \begin{itemize}
\item $df$ maps the  subbundle ${\cal F}_q$ of $TM$ to the subbundle ${\cal
F}_q$ of $TM^\prime$;
 \item the induced map $df: TM/{\cal F}_q\rightarrow  TM^{\prime}/{\cal F}_q$
commutes with the action of the maximal torus $T$ corresponding to
 the point $q$.
  \end{itemize}
\end{definition}
\begin{example}
 Let $M$ be a hypercomplex manifold equipped with a tri-holomorphic action of a
Lie
 group $H$ for which a moment map exists. Such a moment map  $\mu:M\rightarrow
 {\frak h}^\ast\otimes {\Bbb R}^3$ can be viewed as
 a morphism of two GHC-manifolds: $M$ and  ${\frak h}^\ast\otimes {\Bbb
 R}^3$, where the latter one is equipped with the flat $2$-hypercomplex
structure given
 by the action of $SU(2)$ on the second factor.
\end{example}
The following fact follows directly from the definition of morphisms.
\begin{proposition} A map $f:M\rightarrow M^\prime$ of two regular GHC-manifolds is a
 morphism if and only if there is a  holomorphic bundle map $F$
between the twistor spaces $Z,Z^\prime$ of $M$ and $M^\prime$ such that
$\widetilde{f(m)}=F(\tilde{m})$ where $\tilde{p}$ denotes a section of the
twistor space corresponding to a point $p$ of the manifold.
\end{proposition}
We also have:
\begin{proposition}
Let $M$ and $M^\prime$ be two $k$-hypercomplex manifolds. Then any GHC-morphism
$f$ between
$M$ and $M^\prime$ respects the generalised hypercomplex structure, i.e. $df$
commutes
with the action of $SU(2)$. \label{morphism}
\end{proposition}
\begin{pf} This is a local statement, so we can assume that $M$ and $M^\prime$ are regular.
Let $F$ be the a fibrewise mapping of the twistor spaces of $M$ and
 $M^\prime$ given by the previous proposition. Thus $dF$ at a section is a
 linear mapping between $L\otimes {\Bbb C}^n$
 and $L\otimes {\Bbb C}^m$, where $L$ is the line bundle corresponding to
 the representation $V$. Such a mapping is given by a constant linear map
 from ${\Bbb C}^n$ to ${\Bbb C}^m$ and hence $df$ commutes with the action of
 $SU(2)$.
\end{pf}

The following proposition generates a large number of generalised
hypercomplex structures. It is an analogue of the fact \cite{GT} that
$4$-dimensional hypercomplex manifolds with circle symmetry fiber over
hyper-CR-manifolds.
\begin{proposition} Let $M$ be a regular $k$-hypercomplex manifold of dimension
$(k+\nolinebreak 1)^2$ equipped with a free circle action which respects the
$k$-hypercomplex structure and such that the vector field $V$ generated by
the action satisfies the following two conditions:
\begin{itemize}
\item  for any $q\in {\Bbb C}P^1$ and any $m\in M$, $V_m\not\in ({\cal
F}_q)_m$;
\item for any $m$ in $M$, $SU(2)V_m$ linearly generates  all of $T_mM$.
\end{itemize}
Then $M/S^1$ is a $(k+1)$-hypercomplex manifold and the projection
$M\rightarrow M/S^1$ is a morphism of generalised hypercomplex structures.
\end{proposition}
\begin{pf}
Let $Z$ be the twistor space of $M$. The first condition on $V$
implies that there is a free local action of ${\Bbb C}^\ast$ on fibers of
$Z$ and we can form a complex manifold $Z_{\text red}$ by taking fibrewise
quotients (at least locally). Sections of $Z$ descend to sections of
$Z_{\text red}$ and we have to show that they have correct normal bundle.
This normal bundle $N$ fits in the exact sequence $$ 0\rightarrow {\cal
O}\rightarrow {\cal O}(k)^{k+1}\rightarrow N \rightarrow 0.$$ Thus $N$ is a
vector bundle of rank $k$ and the first Chern class $k(k+1)$. In addition,
the second condition on $V$ implies that ${\cal O}$ does not embed into any
proper subbundle of ${\cal O}(k)^{k+1}$ of the form ${\cal O}(k)^s$.
Therefore any line bundle which is a direct summand of $N$ has the first
Chern class at least $k+1$. It follows that $N\simeq {\cal
O}(k+1)^{k}$.\end{pf}

\section{Symplectic $k$-hypercomplex structures\label{form}}

Let $M$ be a generalised hypercomplex manifold.  Let $V$ be the irreducible representation of
$SU(2)$ such that $T_mM\simeq V\otimes {\Bbb R}^n$ for any $m\in M$. If $V^\cx\simeq S^k\cx^2$, then we denote by $V^{[2]}$ the irreducible representation of $SU(2)$ on a real vector space such that 
$(V^{[2]})^{\Bbb C}\simeq  S^{2k}\cx^2$.

\begin{definition} A $k$-symplectic structure on $M$ is a closed nondegenerate $SU(2)$-invariant $2$-form $\omega$ with values in $V^{[2]}$, i.e. $\omega:\Lambda^2TM\rightarrow V^{[2]}$. \end{definition} 
To explain this notion, recall  the decomposition $T^\cx M\simeq E\otimes H$, where $H$ is the trivial bundle with fiber $S^k\cx^2$. The space of $2$-vectors decomposes as $$\Lambda^2
T M^\cx=\bigl(\Lambda^2E\otimes S^2 H\bigr) \oplus
\bigl(S^2E\otimes \Lambda^2 H\bigr).$$ 
Observe that this expression contains exactly one component
$S^{2k}\cx^2$, which lies in the first term. Therefore a $k$-symplectic form $\omega$ is  equivalent to a nondegenerate 2-form
$\omega_E$ on $E$, compatible with the quaternionic or real structure, such that the canonical map $$\Lambda^2 \bigl(T M^\cx\bigr)\longrightarrow
\Lambda^2E\otimes S^{2k}\cx^2\stackrel{\omega_E}{\longrightarrow} S^{2k}\cx^2 $$
defines a closed $H^{[2]}$-valued $2$-form.
\par
For regular GHC-manifolds we can also give an interpretation in terms of the twistor space $Z$.  We recall that $Z$ fibers over ${\Bbb C}P^1$ and $M$ can
be identified with the space of real sections with normal bundle $L\otimes {\Bbb
C}^n$, where $L$ is a given ample line bundle on ${\Bbb C}P^1$. The $k$-symplectic form is equivalent to  a fibrewise
complex-symplectic form $\omega_Z$ on $Z$ respecting  the real structure of $Z$.
We note that $\omega_Z$ takes values in the line bundle $L^2$, whose space of real sections is $V^{[2]}$. Such a definition of a $k$-symplectic structure is given in \cite{DM}.
\begin{remark} Suppose that the form $\omega_E$ is positive on $(e,\tilde{e})$,
where $e\mapsto \tilde{e}$ is the quaternionic or real structure of $E$. If $M$ is a
symplectic $k$-hypercomplex manifold with $k$ odd, then $M$ admits a canonical
Riemannian metric. Indeed, $S^2 TM^\cx$ splits into the direct sum of
$S^2E\otimes S^2H$ and $\Lambda^2E\otimes \Lambda^2H$. Since, for $k$ odd, the
decomposition of $\Lambda^2 S^k\cx^2$ into irreducibles contains a trivial
representation,  $\Lambda^2H$ has a canonical symplectic $2$-form $\omega_H$,
also compatible with the quaternionic structure. Then $\omega_E\otimes\omega_H$
is a nondegenerate symmetric tensor on $TM^\cx$, compatible with the real
structure and giving a metric on $M$.
\par
If $k$ is even, then it is $S^2H$ which contains a trivial representation,
and so we obtain a canonical symplectic form on  $M$, via the decomposition
$\Lambda^2 TM^\cx=\bigl(\Lambda^2E\otimes S^2 H\bigr) \oplus
\bigl(S^2E\otimes \Lambda^2 H\bigr)$.
\end{remark}

Now suppose that a symplectic GHC-manifold admits a proper and free action of a
Lie group $G$
which respects the GHC-structure and the symplectic form. Suppose, in addition,
that the action is Hamiltonian, i.e. there is a $G$-equivariant map
$$\mu:M\rightarrow V^{[2]}\otimes{\frak g}^\ast$$
having the usual property that its differential evaluated on an element $\rho$
of ${\frak g}$ is equal to the contraction of $\omega$ with the vector field
generated by the action of $\exp(t\rho)$.
\par
We can now construct GHC-manifolds using the symplectic quotient
construction, i.e. defining the reduced manifold as $\mu^{-1}(s)/G$, where $s\in
V^{[2]}\otimes{\frak g}^\ast$ is invariant under $G$-action.  For regular manifolds this is
 equivalent to taking the
complex-symplectic quotient along the fibers of $Z$ to obtain a new twistor
space $Z_{\text red}$.
The generic section of $Z_{\text red}$ which descended from
a section of $Z$ will have the correct normal bundle. Globality is however
hard to come by, as even the linear version of the symplectic quotient works
only generically.  Thus, if $T$ is a trivial isotropic $\sigma$-subbundle of
a symplectic $\sigma$-bundle $E={\cal O}(k)\otimes {\Bbb C}^{2n}$, then it is not
always true that $T^\perp/T$ splits into the sum of line bundles of degree
$k$.

\section{Examples\label{examples}}

The symplectic quotient described in the previous section provides many examples
of generalised hypercomplex structures. Thus we start with the flat
$k$-hypercomplex structure on $V={\Bbb R}^{k+1}\otimes{\Bbb R}^{2n}$. The
twistor space of this GHC-manifold is the total space of the vector bundle ${\cal
 O}(k)\otimes{\Bbb C}^{2n}$. It is naturally a symplectic bundle and we consider
the induced $k$-symplectic structure on $V$. Now we choose a group $H$ preserving
the symplectic GHC-structure. $H$ must be a subgroup of $Sp(n,{\Bbb C})$ and its
representation on ${\Bbb C}^{2n}$ is quaternionic if $k$ is odd and real if $k$
is even. Thus $H$ is a subgroup of $Sp(n)$ if $n$ is odd and of $Sp(2n,{\Bbb
R})$ if $n$ is even. The $k$-symplectic quotient $M$ of $V$ by $H$  will have a
(symplectic) $k$-hypercomplex manifold on an open dense subset. We observe that
the fibers of the twistor space of $M$ look generically the same as for the
twistor space of the hypercomplex manifold obtained by the same quotient
construction with $k$ replaced by $1$. Thus we can think of resulting manifolds as  $k$-hypercomplex analogues
of ALE-spaces, complex coadjoint orbits, toric hyperk\"ahler manifolds etc.
\par
We shall consider two  examples in greater detail.
\subsection{$k$-Eguchi-Hanson manifold}

We begin with $\cx^{k+1}\otimes \cx^4$ which we view as a complexified  flat $k$-hypercomplex manifold. We have the ${\Bbb C}^\ast$-action on the second factor $t\cdot(z_1,z_2,w_1,w_2)= (tz_1,tz_2,t^{-1}w_1,t^{-1}w_2)$ which induces an action of $S^1$ on the underlying $k$-hypercomplex space $V$ which is ${\Bbb R}^{2k+2}\otimes {\Bbb R}^2$ for odd $k$, and  ${\Bbb R}^{k+1}\otimes {\Bbb R}^4$ for even $k$. We consider the $k$-symplectic quotient of this space by the $S^1$-action. According to the discussion in the previous section, this is equivalent to taking the fibrewise symplectic quotient of ${\cal O}(k)\otimes {\Bbb C}^4$ by  ${\Bbb C}^\ast$. In other words, we identify $V^{\cx}$ with the space of sections of  ${\cal O}(k)\otimes {\Bbb C}^4$, i.e. with quadruples of polynomials $\bigl(z_1(\zeta),z_2(\zeta),w_1(\zeta),w_2(\zeta)\bigr)$ of degree $k$. The moment map $\mu: V^{\cx}\rightarrow H^0({\cal O}(2k))$ is simply 
$$\mu\bigl(z_1(\zeta),z_2(\zeta),w_1(\zeta),w_2(\zeta)\bigr)= z_1(\zeta)w_1(\zeta)+ z_2(\zeta)w_2(\zeta).$$
It is clear that a  section $s$ of ${\cal O}(2k)$ will be a regular value for $\mu$ as soon as $s$ does not have double zeros. This way we obtain a manifold $\mu^{-1}(s)/{\Bbb C}^\ast$. If we begin with real sections of ${\cal O}(k)\otimes {\Bbb C}^4$ and a real section $s$ and quotient by $S^1$, we shall obtain a manifold $M_k$ of dimension $2k+2$. This manifold carries a $k$-hypercomplex structure on an open dense subset and one can call it the {\em $k$-Eguchi-Hanson manifold}. 
\par 
We shall describe $M_{k-1}$ or rather its complexification  $M_{k-1}^\cx$. Since the group (circle) is abelian, the moment map is invariant and we can take first the GIT-quotient of $\cx^{4k}$ by $\cx^\ast$ and then its intersection with a level set of the moment map, which becomes linear. In the present case the GIT quotient of $\cx^{4k}$ by $\cx^\ast$ is simply the variety of  $2k\times 2k$-matrices of rank $1$. Let us choose coordinates so that the polynomials $z_p(\zeta),w_p(\zeta)$,  $p=1,2$, are of the form:
$$ z_1(\zeta)=\sum_{i=1}^{k} z_{i}\zeta^{i-1}, \qquad  z_2(\zeta)=\sum_{i=k+1}^{2k+2} z_{i}\zeta^{i-1}, $$
$$w_1(\zeta)=\sum_{i=1}^{k} w_i\zeta^{k-i},\qquad w_2(\zeta)=\sum_{i=1}^{k} w_{k+i}\zeta^{2k-i}.$$
Then $M_{k-1}^\cx$ is the intersection of the variety of rank $1$ matrices $[a_{ij}]$ with the affine subspace described by equations
$$\sum_{i=1}^{k-p} a_{i,p+i}+ \sum_{i=1}^{k-p} a_{k+i,k+p+i}
=\tau_i, \enskip p=0,\dots,k-1,$$
$$\sum_{i=1}^{k-p} a_{p+i,i}+ \sum_{i=1}^{k-p} a_{k+p+i,k+i}
=\nu_i, \enskip p=0,\dots,k-1.$$
These equations fix the sum of entries of the matrix lying in the diagonal $(k\times k)$-blocks and parallel to the main diagonal. To obtain $M$ we impose the reality condition $\bar{w}_i=z_i$. The $k$-hypercomplex structure on $M$ should be equivalent to the one obtained in \cite{DM2} via a solution to the Pleba\'nski heavenly equation.

\subsection{Infinite-dimensional quotients}

Just as in case of hypercomplex manifolds, we can consider infinite-dimensional $k$-symplectic quotients. We shall discuss the case of $k=2$. Let $G$ be a compact semisimple Lie group with Lie algebra $\g$. Consider the space ${\cal A}$ of $\g$-valued smooth functions $T_0,\dots T_5$ on an interval $[0,1]$. These can be viewed as real sections of the bundle $C^\infty\bigl([0,1], \g^\cx\bigr)\otimes ({\cal O}(2)\oplus {\cal O}(2))$ on $\cx P^1$, and so we have a flat infinite-dimensional symplectic $2$-hypercomplex manifold. The gauge group of transformations $g:[0,1]\rightarrow \g$ acts via:
$$T_0\mapsto \ad (g) T_0 -\dot{g}g^1,\quad T_i\mapsto \ad (g) T_i,\enskip i=1,\dots,5.$$
The $2$-symplectic quotient by the gauge subgroup of $g$ with $g(0)=g(1)=1$ will be described by an analogue of Nahm's equations. The simplest way to describe the equations is again to perform the quotient along the fibers of the twistor space. Let us write 
$$\beta=T_4+iT_5,\quad \gamma=T_2+iT_3, \quad\alpha=T_0+iT_1,$$  
and put
$$B(\zeta)=\beta +\gamma \zeta +(\alpha+\alpha^\ast)\zeta^2 -\gamma^\ast \zeta^3 +\beta^\ast\zeta^4,$$
$$ A(\zeta)=\alpha -\gamma^\ast \zeta +\beta^\ast\zeta^2.$$
Here $A$ and $-A+B/\zeta^2$ are sections of $C^\infty\bigl([0,1],\g^\cx\bigr) \otimes {\cal O}(2)$ and provide coordinates on  ${\cal A}$. Alternatively, we can choose $B$ and $T_0$ as coordinates on ${\cal A}$, which corresponds to writing ${\cal O}(2)\oplus {\cal O}(2)$ as the extension
$$0\rightarrow {\cal O}\rightarrow {\cal O}(2)\oplus {\cal O}(2)\rightarrow {\cal O}(4)\rightarrow 0.$$
The $2$-symplectic quotient $M$ is then described as the moduli space of solutions to
\begin{equation}\dot{B}(\zeta)=[B(\zeta),A(\zeta)],\label{N0} \end{equation}
which is simply the space of sections of the fibrewise symplectic quotient of the twistor space $C^\infty\bigl([0,1], \g^\cx\bigr)\otimes ({\cal O}(2)\oplus {\cal O}(2))$ . 
\par
The equations \eqref{N0} are equivalent to the following:
\begin{equation}\frac{d}{dt}{\beta}=[\beta,\alpha],\label{N1}\end{equation}
\begin{equation}\frac{d}{dt}{\gamma}=[\gamma,\alpha] -[\beta,\gamma^\ast],\label{N2}\end{equation}
\begin{equation}\frac{d}{dt}(\alpha+\alpha^\ast)=[\alpha^\ast,\alpha]+[\beta,\beta^\ast]-[\gamma,\gamma^\ast].\label{N3}\end{equation}
The gauge group acts via:
$$\alpha\mapsto \ad (g)\alpha-\dot{g}g^{-1},\quad \beta\mapsto \ad (g)\beta,\quad \gamma\mapsto \ad(g)\gamma.$$
To identify $M$, use a gauge transformation $g$, with $g(0)=1$, to make $\alpha$ hermitian. Then the map $$[\alpha,\beta,\gamma]\longmapsto \bigl(g(1),\alpha(0)+\alpha^\ast(0),\beta(0),\gamma(0)\bigr)$$
is diffeomorphism from $M$ into $G\times \g\times\g^\cx\times\g^\cx$. We observe that $M$ fibers over  a complex manifold $Z_0$ consisting of solutions to $\dot{\beta}=[\beta,\alpha]$ modulo {\em complex} gauge transformations. This space is easily identified  with $G^\cx\times g^\cx\simeq T^\ast G^\cx$ (one can make $\alpha= 0$ via a complex gauge transformation $g$ with $g(0)=1$). 
\par
The $2$-hypercomplex structure is $G$-invariant and symplectic.
It would be interesting to identify the hypercomplex extension of $M$, given by Theorem \ref{EXT}. An obvious candidate is the hypercomplex structure on a neighbourhood of the zero section of $T^\ast H^\cx$, where $H$ is the tangent group of $G$, i.e. the semidirect product of $G$ and $\g$.
\par
Even more interesting would be to consider the $3$-hypercomplex manifold given by a similar construction and identify the hypercomplex structures, described by Theorem \ref{k-2}, on real $\alpha$-surfaces.
\par
Finally, one can impose different boundary conditions (poles at $t=0$, conjugacy classes at infinity, etc.) on equations \eqref{N1}-\eqref{N3} and their $k$-hypercomplex, $k>2$, analogues. This yields not only many $k$-hypercomplex structures, but also, in view of Theorems \ref{k-2} and \ref{EXT}, many hypercomplex structures.

\medskip

{\bf Acknowledgment}. I thank David Calderbank and Misha  Verbitsky for reading the manuscript and for comments.

\end{document}